\documentclass[11pt]{amsart} 
\usepackage{amsmath,amssymb,amsthm}
\newtheorem{theorem}{Theorem}
\newtheorem{proposition}[theorem]{Proposition}

\newtheorem{lemma}[theorem]{Lemma}

\begin{document}

\title{Embedded minimal tori in $S^3$ and the Lawson conjecture} 
\author{Simon Brendle}
\begin{abstract}
We show that any embedded minimal torus in $S^3$ is congruent to the Clifford torus. This answers a question posed by H.B.~Lawson, Jr., in 1970.
\end{abstract}
\address{Department of Mathematics \\ Stanford University \\ Stanford, CA 94305}
\thanks{The author was supported in part by the National Science Foundation under grants DMS-0905628 and DMS-1201924.}
\maketitle

\section{Introduction}

The study of minimal surfaces is one of the oldest subjects in differential geometry. Of particular interest are minimal surfaces in spaces of constant curvature, such as the Euclidean space $\mathbb{R}^3$ or the sphere $S^3$.  The case of the sphere $S^3$ turns out to be very interesting: for example, while there are no closed minimal surfaces in $\mathbb{R}^3$, the sphere $S^3$ does contain closed minimal surfaces. The simplest example of a minimal surface in $S^3$ is the equator. Another basic example is the so-called Clifford torus. Identifying $S^3$ with the unit sphere in $\mathbb{R}^4$, the Clifford torus is defined by 
\[\Big \{ (x_1,x_2,x_3,x_4) \in S^3: x_1^2+x_2^2=x_3^2+x_4^2=\frac{1}{2} \Big \}.\] 
We note that the principal curvatures of the Clifford torus are $1$ and $-1$, and the intrinsic Gaussian curvature vanishes identically. 

In 1970, Lawson \cite{Lawson2} proved that, given any positive integer $g$, there exists at least one compact embedded minimal surface in $S^3$ with genus $g$ (cf. \cite{Lawson2}, Section 6). Moreover, he showed that there are at least two such surfaces unless the genus $g$ is a prime number. Additional examples of compact embedded minimal surfaces in $S^3$ were later found by Karcher, Pinkall, and Sterling \cite{Karcher-Pinkall-Sterling} and, more recently, by Kapouleas and Yang \cite{Kapouleas-Yang}. The construction of Karcher, Pinkall, and Sterling uses tesselations of $S^3$ into cells that have the symmetry of a Platonic solid in $\mathbb{R}^3$; the resulting minimal surfaces have genus $3$, $5$, $6$, $7$, $11$, $19$, $73$, and $601$, respectively. The result of Kapouleas and Yang relies on a so-called doubling construction: roughly speaking, this construction involves joining together two nearby copies of the Clifford torus by a large number of catenoid necks. The resulting surfaces have small mean curvature, and Kapouleas and Yang employed the implicit function theorem to deform these surfaces to exact solutions of the minimal surface equation. We note that Kapouleas has recently described a similar doubling construction involving the equator instead of the Clifford torus (cf. \cite{Kapouleas-survey}, Section 2.4).

It was shown by Almgren in 1966 that any immersed minimal two-sphere in $S^3$ is totally geodesic, and therefore congruent to the equator (see \cite{Almgren}, p.~279). Almgren's proof uses the method of Hopf differentials (see also \cite{Lawson2}, Proposition 1.5). In 1970, Lawson \cite{Lawson3} conjectured that the Clifford torus is the only compact embedded minimal surface in $S^3$ of genus $1$. In this paper, we give an affirmative answer to Lawson's conjecture:

\begin{theorem}
\label{main.theorem}
Suppose that $F: \Sigma \to S^3$ is an embedded minimal torus in $S^3$. Then $F$ is congruent to the Clifford torus.
\end{theorem}

We note that the embeddedness assumption in Theorem \ref{main.theorem} is crucial: in fact, Lawson \cite{Lawson2} has constructed an infinite family of minimal immersions from the torus and the Klein bottle into $S^3$ (see also \cite{Hsiang-Lawson}).

Lawson's conjecture has attracted considerable interest over the past decades, and various partial results are known. For example, it was shown by Urbano \cite{Urbano} that any minimal torus in $S^3$ which has Morse index at most $5$ is congruent to the Clifford torus. Moreover, Ros \cite{Ros} was able to verify Lawson's conjecture for surfaces that are invariant under reflection across each coordinate plane. Montiel and Ros \cite{Montiel-Ros} linked Lawson's conjecture to a conjecture of Yau (cf. \cite{Yau}) concerning the first eigenvalue of the Laplacian on a minimal surface. Yau's conjecture is discussed in more detail in \cite{Choe}, \cite{Choe-Soret}, and \cite{Choi-Wang}. Finally, we note that Marques and Neves recently showed that the Clifford torus has smallest area among all minimal surfaces in $S^3$ of genus at least $1$ (cf. \cite{Marques-Neves}, Theorem B). The method used in \cite{Marques-Neves} is completely different from ours; it relies on the min-max theory for minimal surfaces and the rigidity theorem of Urbano.

Our method of proof is inspired in part by the pioneering work of G.~Huisken \cite{Huisken} on the curve shortening flow, as well as by recent work of B.~Andrews \cite{Andrews} on the mean curvature flow. Let us digress briefly to review these results.

Given a one-parameter family of embedded curves $F_t: S^1 \to \mathbb{R}^2$, Huisken considered the quantity 
\[W_t(x,y) = \frac{L(t)}{|F_t(x) - F_t(y)|} \, \sin \Big ( \frac{\pi \, d_t(x,y)}{L(t)} \Big ),\] 
where $L(t)$ denotes the total length of the curve $F_t$ and $d_t(x,y)$ denotes the intrinsic distance of two points $x,y \in S^1$. Huisken discovered that, if the curves $F_t$ evolve by the curve shortening flow, then the supremum of the function $W_t(x,y)$ is monotone decreasing in $t$. Using this monotonicity formula, Huisken was able to show that the curve shortening flow deforms any embedded curve in the plane to a round point, thereby giving a direct proof of a theorem of Grayson \cite{Grayson}.

Huisken's technique was developed further in a recent paper by B.~Andrews \cite{Andrews}. Andrews considered a one-parameter family of embedded hypersurfaces $F_t: M \to \mathbb{R}^{n+1}$ which have positive mean curvature and evolve by the mean curvature flow. By applying the maximum principle to a suitable function $W_t(x,y)$ defined on $M \times M$, Andrews obtained a new proof of the noncollapsing estimate established earlier by Sheng and Wang \cite{Sheng-Wang} (see also \cite{White}). The argument in \cite{Andrews} relies in a crucial way on the positivity of the mean curvature; in particular, the argument is not applicable in the case of minimal surfaces. 

We now describe the main ideas involved in the proof of Theorem \ref{main.theorem}. Let $F: \Sigma \to S^3$ be an embedded minimal torus in $S^3$. It follows from work of Lawson that the surface $\Sigma$ has no umbilic points. Consequently, the quantity 
\[\kappa = \sup_{x,y \in \Sigma, \, x \neq y} \sqrt{2} \, \frac{|\langle \nu(x),F(y) \rangle|}{|A(x)| \, (1 - \langle F(x),F(y) \rangle)}\] 
is finite. If $\kappa \leq 1$, we can show that the second fundamental form of $F$ is parallel. From this, we deduce that the induced metric on $\Sigma$ is flat. A classical theorem of Lawson \cite{Lawson1} then implies that $F$ is congruent to the Clifford torus.

Hence, it remains to consider the case $\kappa > 1$. In order to handle this case, we apply the maximum principle to the function 
\begin{equation} 
Z(x,y) = \frac{\kappa}{\sqrt{2}} \, |A(x)| \, (1 - \langle F(x),F(y) \rangle) + \langle \nu(x),F(y) \rangle. 
\end{equation}
By definition of $\kappa$, the function $Z(x,y)$ is nonnegative for all points $x,y \in \Sigma$. Moreover, after replacing $\nu$ by $-\nu$ if necessary, we can find two points $\bar{x},\bar{y} \in \Sigma$ such that $\bar{x} \neq \bar{y}$ and $Z(\bar{x},\bar{y}) = 0$. Since the function $Z$ attains its global minimum at $(\bar{x},\bar{y})$, the first derivatives of the function $Z$ at the point $(\bar{x},\bar{y})$ vanish, and the Hessian of the function $Z$ at the point $(\bar{x},\bar{y})$ is nonnegative. In order to analyze the Hessian of the function $Z$, we use the identity 
\begin{equation} 
\label{simons.id}
\Delta_\Sigma(|A|) - \frac{\big | \nabla |A| \big |^2}{|A|} + (|A|^2 - 2) \, |A| = 0.
\end{equation} 
The relation (\ref{simons.id}) is a consequence of the classical Simons identity (cf. \cite{Simons}).

At this point, we encounter a major obstacle: the identity (\ref{simons.id}) contains a gradient term which has an unfavorable sign. However, by exploiting special identities arising from the first variation of the function $Z(x,y)$, we are able to extract a gradient term which has a favorable sign (see Proposition \ref{key.ingredient} below). Surprisingly, this term precisely offsets the bad term coming from the Simons identity! It is this insight which makes the maximum principle work. This leads to the inequality 
\begin{equation} 
\label{a}
\sum_{i=1}^2 \frac{\partial^2 Z}{\partial x_i^2}(\bar{x},\bar{y}) \leq \sqrt{2} \, \kappa \, |A(\bar{x})|. 
\end{equation}
Moreover, we compute 
\begin{equation} 
\label{b}
\sum_{i=1}^2 \frac{\partial^2 Z}{\partial y_i^2}(\bar{x},\bar{y}) = \sqrt{2} \, \kappa \, |A(\bar{x})|. 
\end{equation}
In order to absorb the terms $\sqrt{2} \, \kappa \, |A(\bar{x})|$ on the right hand side of (\ref{a}) and (\ref{b}), we consider the mixed partial derivatives $\frac{\partial^2 Z}{\partial x_i \, \partial y_i}$. It turns out that 
\begin{equation} 
\label{c}
\sum_{i=1}^2 \frac{\partial^2 Z}{\partial x_i \, \partial y_i} = -\sqrt{2} \, \kappa \, |A(\bar{x})|, 
\end{equation}
where $(x_1,x_2)$ and $(y_1,y_2)$ are suitably chosen local coordinates around $\bar{x}$ and $\bar{y}$, respectively. By combining (\ref{a}), (\ref{b}), and (\ref{c}), we can make the terms $\sqrt{2} \, \kappa \, |A(\bar{x})|$ cancel, and we obtain 
\begin{equation} 
\label{main.inequality}
\sum_{i=1}^2 \frac{\partial^2 Z}{\partial x_i^2}(\bar{x},\bar{y}) + 2 \sum_{i=1}^2 \frac{\partial^2 Z}{\partial x_i \, \partial y_i}(\bar{x},\bar{y}) + \sum_{i=1}^2 \frac{\partial^2 Z}{\partial y_i^2}(\bar{x},\bar{y}) \leq 0. 
\end{equation}
We now apply the strict maximum principle for degenerate elliptic equations (cf. \cite{Bony}) to the function $Z(x,y)$. From this, we deduce that the function $|A|$ is constant. This again implies that $F$ is congruent to the Clifford torus.

The idea of exploiting the mixed partial derivatives $\frac{\partial^2 Z}{\partial x_i \, \partial y_i}$ goes back to work of Huisken \cite{Huisken} and was also used in \cite{Andrews}. An interesting feature of our argument is that we need to use the full strength of the mixed partial derivative terms $\frac{\partial^2 Z}{\partial x_i \, \partial y_i}$.

It is a pleasure to thank Professors Gerhard Huisken and Brian White for discussions. 

\section{The key technical ingredient}

\label{calc}

Let $F: \Sigma \to S^3$ be an embedded minimal surface in $S^3$ (viewed as the unit sphere in $\mathbb{R}^4$). Moreover, let $\Phi$ be a positive function on $\Sigma$. We consider the expression 
\[Z(x,y) = \Phi(x) \, (1 - \langle F(x),F(y) \rangle) + \langle \nu(x),F(y) \rangle.\] 
Let us consider a pair of points $\bar{x} \neq \bar{y}$ with the property that $Z(\bar{x},\bar{y}) = 0$ and the differential of $Z$ at the point $(\bar{x},\bar{y})$ vanishes. Let $(x_1,x_2)$ be geodesic normal coordinates around $\bar{x}$, and let $(y_1,y_2)$ be geodesic normal coordinates around $\bar{y}$. 

At the point $(\bar{x},\bar{y})$, we have 
\begin{align} 
\label{first.order.condition.a}
0 = \frac{\partial Z}{\partial x_i}(\bar{x},\bar{y}) 
&= \frac{\partial \Phi}{\partial x_i}(\bar{x}) \, (1 - \langle F(\bar{x}),F(\bar{y}) \rangle) \notag \\ 
&- \Phi(\bar{x}) \, \Big \langle \frac{\partial F}{\partial x_i}(\bar{x}),F(\bar{y}) \Big \rangle + h_i^k(\bar{x}) \, \Big \langle \frac{\partial F}{\partial x_k}(\bar{x}),F(\bar{y}) \Big \rangle 
\end{align}
and 
\begin{equation} 
\label{first.order.condition.b}
0 = \frac{\partial Z}{\partial y_i}(\bar{x},\bar{y}) = -\Phi(\bar{x}) \, \Big \langle F(\bar{x}),\frac{\partial F}{\partial y_i}(\bar{y}) \Big \rangle + \Big \langle \nu(\bar{x}),\frac{\partial F}{\partial y_i}(\bar{y}) \Big \rangle. 
\end{equation}
We will make extensive use of these relations in the subsequent arguments.

Without loss of generality, we may assume that the second fundamental form at $\bar{x}$ is diagonal, so that $h_{11}(\bar{x}) = \lambda_1$, $h_{12}(\bar{x}) = 0$, and $h_{22}(\bar{x}) = \lambda_2$. We denote by $w_i$ the reflection of the vector $\frac{\partial F}{\partial x_i}(\bar{x})$ across the hyperplane orthogonal to $F(\bar{x})-F(\bar{y})$, i.e. 
\[w_i = \frac{\partial F}{\partial x_i}(\bar{x}) - 2 \, \Big \langle \frac{\partial F}{\partial x_i}(\bar{x}),\frac{F(\bar{x})-F(\bar{y})}{|F(\bar{x})-F(\bar{y})|} \Big \rangle \, \frac{F(\bar{x})-F(\bar{y})}{|F(\bar{x})-F(\bar{y})|}.\] 
By a suitable choice of the coordinate system $(y_1,y_2)$, we can arrange that $\big \langle w_1,\frac{\partial F}{\partial y_1}(\bar{y}) \big \rangle \geq 0$, $\big \langle w_1,\frac{\partial F}{\partial y_2}(\bar{y}) \big \rangle = 0$, and $\big \langle w_2,\frac{\partial F}{\partial y_2}(\bar{y}) \big \rangle \geq 0$. 

\begin{lemma}
The vectors $F(\bar{y})$ and $\Phi(\bar{x}) \, F(\bar{x}) - \nu(\bar{x})$ are linearly independent. 
\end{lemma}

\textbf{Proof.} 
Using the identity 
\[\langle \Phi(\bar{x}) \, F(\bar{x}) - \nu(\bar{x}),F(\bar{y}) \rangle = \Phi(\bar{x}) - Z(\bar{x},\bar{y}) = \Phi(\bar{x}),\] 
we obtain 
\begin{align*} 
&|\Phi(\bar{x}) \, F(\bar{x}) - \nu(\bar{x})|^2 \, |F(\bar{y})|^2 - \langle \Phi(\bar{x}) \, F(\bar{x}) - \nu(\bar{x}),F(\bar{y}) \rangle^2 \\ 
&= |\Phi(\bar{x}) \, F(\bar{x}) - \nu(\bar{x})|^2 - \Phi(\bar{x})^2 = 1. 
\end{align*} From this, the assertion follows. \\

\begin{lemma}
\label{algebra}
We have $w_1 = \frac{\partial F}{\partial y_1}(\bar{y})$ and $w_2 = \frac{\partial F}{\partial y_2}(\bar{y})$.
\end{lemma} 

\textbf{Proof.} 
A straightforward calculation gives 
\begin{align*} 
\langle w_i,F(\bar{y}) \rangle 
&= \Big \langle \frac{\partial F}{\partial x_i}(\bar{x}),F(\bar{y}) \Big \rangle \\ 
&+ 2 \, \Big \langle \frac{\partial F}{\partial x_i}(\bar{x}),F(\bar{y}) \Big \rangle \, \frac{\langle F(\bar{x})-F(\bar{y}),F(\bar{y}) \rangle}{|F(\bar{x})-F(\bar{y})|^2} = 0 
\end{align*}
and 
\begin{align*} 
&\langle w_i,\Phi(\bar{x}) \, F(\bar{x}) - \nu(\bar{x}) \rangle \\ 
&= 2 \, \Big \langle \frac{\partial F}{\partial x_i}(\bar{x}),F(\bar{y}) \Big \rangle \, \frac{\langle F(\bar{x})-F(\bar{y}),\Phi(\bar{x}) \, F(\bar{x}) - \nu(\bar{x}) \rangle}{|F(\bar{x})-F(\bar{y})|^2} \\ 
&= 2 \, \Big \langle \frac{\partial F}{\partial x_i}(\bar{x}),F(\bar{y}) \Big \rangle \, \frac{Z(\bar{x},\bar{y})}{|F(\bar{x})-F(\bar{y})|^2} \\ 
&= 0. 
\end{align*}
On the other hand, the vectors $\frac{\partial F}{\partial y_1}(\bar{y})$ and $\frac{\partial F}{\partial y_2}(\bar{y})$ satisfy 
\[\Big \langle \frac{\partial F}{\partial y_i}(\bar{y}),F(\bar{y}) \Big \rangle = 0\] 
and 
\[\Big \langle \frac{\partial F}{\partial y_i}(\bar{y}),\Phi(\bar{x}) \, F(\bar{x}) - \nu(\bar{x}) \Big \rangle = -\frac{\partial Z}{\partial y_i}(\bar{x},\bar{y}) = 0.\] 
Since the vectors $F(\bar{y})$ and $\Phi(\bar{x}) \, F(\bar{x}) - \nu(\bar{x})$ are linearly independent, we conclude that the plane spanned by $w_1$ and $w_2$ coincides with the plane spanned by $\frac{\partial F}{\partial y_1}(\bar{y})$ and $\frac{\partial F}{\partial y_2}(\bar{y})$. Moreover, $w_1$ and $w_2$ are orthonormal. Since $\big \langle w_1,\frac{\partial F}{\partial y_2}(\bar{y}) \big \rangle = 0$, we conclude that $w_1 = \pm \frac{\partial F}{\partial y_1}(\bar{y})$ and $w_2 = \pm \frac{\partial F}{\partial y_2}(\bar{y})$. Since $\big \langle w_1,\frac{\partial F}{\partial y_1}(\bar{y}) \big \rangle \geq 0$ and $\big \langle w_2,\frac{\partial F}{\partial y_2}(\bar{y}) \big \rangle \geq 0$, the assertion follows. \\

We next consider the second order derivatives of $Z$ at the point $(\bar{x},\bar{y})$. 

\begin{proposition} 
\label{second.derivatives.1}
We have 
\begin{align*} 
&\sum_{i=1}^2 \frac{\partial^2 Z}{\partial x_i^2}(\bar{x},\bar{y}) \\ 
&= \Big ( \Delta_\Sigma \Phi(\bar{x}) - \frac{|\nabla \Phi(\bar{x})|^2}{\Phi(\bar{x})} + (|A(\bar{x})|^2-2) \, \Phi(\bar{x}) \Big ) \, (1 - \langle F(\bar{x}),F(\bar{y}) \rangle) + 2 \, \Phi(\bar{x}) \\ 
&- \frac{2 \, \Phi(\bar{x})^2 - |A(\bar{x})|^2}{2 \, \Phi(\bar{x}) \, (1-\langle F(\bar{x}),F(\bar{y}) \rangle)} \sum_{i=1}^2 \Big \langle \frac{\partial F}{\partial x_i}(\bar{x}),F(\bar{y}) \Big \rangle^2. 
\end{align*} 
\end{proposition}

\textbf{Proof.} 
It follows from the Codazzi equations that 
\[\sum_{i=1}^2 \frac{\partial}{\partial x_i} h_i^k(\bar{x}) = 0.\] 
This implies 
\begin{align*} 
&\sum_{i=1}^2 \frac{\partial^2 Z}{\partial x_i^2}(\bar{x},\bar{y}) \\ 
&= \sum_{i=1}^2 \frac{\partial^2 \Phi}{\partial x_i^2}(\bar{x}) \, (1 - \langle F(\bar{x}),F(\bar{y}) \rangle) - 2 \sum_{i=1}^2 \frac{\partial \Phi}{\partial x_i}(\bar{x}) \, \Big \langle \frac{\partial F}{\partial x_i}(\bar{x}),F(\bar{y}) \Big \rangle \\ 
&+ 2 \, \Phi(\bar{x}) \, \langle F(\bar{x}),F(\bar{y}) \rangle - |A(\bar{x})|^2 \, \langle \nu(\bar{x}),F(\bar{y}) \rangle \\ 
&= \big ( \Delta_\Sigma \Phi(\bar{x}) + (|A(\bar{x})|^2-2) \, \Phi(\bar{x}) \big ) \, (1 - \langle F(\bar{x}),F(\bar{y}) \rangle) + 2 \, \Phi(\bar{x}) \\ 
&- 2 \sum_{i=1}^2 \frac{\partial \Phi}{\partial x_i}(\bar{x}) \, \Big \langle \frac{\partial F}{\partial x_i}(\bar{x}),F(\bar{y}) \Big \rangle.  
\end{align*} 
Rearranging terms gives 
\begin{align*} 
&\sum_{i=1}^2 \frac{\partial^2 Z}{\partial x_i^2}(\bar{x},\bar{y}) \\ 
&= \Big ( \Delta_\Sigma \Phi(\bar{x}) - \frac{|\nabla \Phi(\bar{x})|^2}{\Phi(\bar{x})} + (|A(\bar{x})|^2-2) \, \Phi(\bar{x}) \Big ) \, (1 - \langle F(\bar{x}),F(\bar{y}) \rangle) + 2 \, \Phi(\bar{x}) \\ 
&+ \frac{1}{\Phi(\bar{x}) \, (1-\langle F(\bar{x}),F(\bar{y}) \rangle)} \\ 
&\hspace{15mm} \cdot \sum_{i=1}^2 \bigg ( \frac{\partial \Phi}{\partial x_i}(\bar{x}) \, (1 - \langle F(\bar{x}),F(\bar{y}) \rangle) - \Phi(\bar{x}) \, \Big \langle \frac{\partial F}{\partial x_i}(\bar{x}),F(\bar{y}) \Big \rangle \bigg )^2 \\ 
&- \frac{\Phi(\bar{x})}{1-\langle F(\bar{x}),F(\bar{y}) \rangle} \sum_{i=1}^2 \Big \langle \frac{\partial F}{\partial x_i}(\bar{x}),F(\bar{y}) \Big \rangle^2. 
\end{align*} 
Using the identity (\ref{first.order.condition.a}), we obtain 
\begin{align*} 
&\sum_{i=1}^2 \frac{\partial^2 Z}{\partial x_i^2}(\bar{x},\bar{y}) \\ 
&= \Big ( \Delta_\Sigma \Phi(\bar{x}) - \frac{|\nabla \Phi(\bar{x})|^2}{\Phi(\bar{x})} + (|A(\bar{x})|^2-2) \, \Phi(\bar{x}) \Big ) \, (1 - \langle F(\bar{x}),F(\bar{y}) \rangle) + 2 \, \Phi(\bar{x}) \\ 
&+ \frac{1}{\Phi(\bar{x}) \, (1-\langle F(\bar{x}),F(\bar{y}) \rangle)} \sum_{i=1}^2 \lambda_i^2 \, \Big \langle \frac{\partial F}{\partial x_i}(\bar{x}),F(\bar{y}) \Big \rangle^2 \\ 
&- \frac{\Phi(\bar{x})}{1-\langle F(\bar{x}),F(\bar{y}) \rangle} \sum_{i=1}^2 \Big \langle \frac{\partial F}{\partial x_i}(\bar{x}),F(\bar{y}) \Big \rangle^2. 
\end{align*} 
Since $\lambda_1^2 = \lambda_2^2 = \frac{1}{2} \, |A(\bar{x})|^2$, the assertion follows. \\

\begin{proposition}
\label{second.derivatives.2}
We have 
\[\frac{\partial^2 Z}{\partial x_i \, \partial y_i}(\bar{x},\bar{y}) = \lambda_i - \Phi(\bar{x}).\] 
\end{proposition}

\textbf{Proof.} 
Using the relation (\ref{first.order.condition.a}) and Lemma \ref{algebra}, we obtain 
\begin{align*} 
&\frac{\partial^2 Z}{\partial x_i \, \partial y_i}(\bar{x},\bar{y}) \\ 
&= -\frac{\partial \Phi}{\partial x_i}(\bar{x}) \, \Big \langle F(\bar{x}),\frac{\partial F}{\partial y_i}(\bar{y}) \Big \rangle + (\lambda_i - \Phi(\bar{x})) \, \Big \langle \frac{\partial F}{\partial x_i}(\bar{x}),\frac{\partial F}{\partial y_i}(\bar{y}) \Big \rangle \\ 
&= \frac{1}{1 - \langle F(\bar{x}),F(\bar{y}) \rangle} \, (\lambda_i - \Phi(\bar{x})) \, \Big \langle \frac{\partial F}{\partial x_i}(\bar{x}),F(\bar{y}) \Big \rangle \, \Big \langle F(\bar{x}),\frac{\partial F}{\partial y_i}(\bar{y}) \Big \rangle \\ 
&+ (\lambda_i - \Phi(\bar{x})) \, \Big \langle \frac{\partial F}{\partial x_i}(\bar{x}),\frac{\partial F}{\partial y_i}(\bar{y}) \Big \rangle \\ 
&= -2 \, (\lambda_i - \Phi(\bar{x})) \, \Big \langle \frac{\partial F}{\partial x_i}(\bar{x}),\frac{F(\bar{x})-F(\bar{y})}{|F(\bar{x})-F(\bar{y})|} \Big \rangle \, \Big \langle \frac{F(\bar{x})-F(\bar{y})}{|F(\bar{x})-F(\bar{y})|},\frac{\partial F}{\partial y_i}(\bar{y}) \Big \rangle \\ 
&+ (\lambda_i - \Phi(\bar{x})) \, \Big \langle \frac{\partial F}{\partial x_i}(\bar{x}),\frac{\partial F}{\partial y_i}(\bar{y}) \Big \rangle \\ 
&= (\lambda_i - \Phi(\bar{x})) \, \Big \langle w_i,\frac{\partial F}{\partial y_i}(\bar{y}) \Big \rangle \\ 
&= \lambda_i - \Phi(\bar{x}), 
\end{align*} 
as claimed. \\

\begin{proposition} 
\label{key.ingredient}
We have
\begin{align*} 
&\sum_{i=1}^2 \frac{\partial^2 Z}{\partial x_i^2}(\bar{x},\bar{y}) + 2 \sum_{i=1}^2 \frac{\partial^2 Z}{\partial x_i \, \partial y_i}(\bar{x},\bar{y}) + \sum_{i=1}^2 \frac{\partial^2 Z}{\partial y_i^2}(\bar{x},\bar{y}) \\ 
&= \Big ( \Delta_\Sigma \Phi(\bar{x}) - \frac{|\nabla \Phi(\bar{x})|^2}{\Phi(\bar{x})} + (|A(\bar{x})|^2-2) \, \Phi(\bar{x}) \Big ) \, (1 - \langle F(\bar{x}),F(\bar{y}) \rangle) \\ 
&- \frac{2 \, \Phi(\bar{x})^2 - |A(\bar{x})|^2}{2 \, \Phi(\bar{x}) \, (1-\langle F(\bar{x}),F(\bar{y}) \rangle)} \sum_{i=1}^2 \Big \langle \frac{\partial F}{\partial x_i}(\bar{x}),F(\bar{y}) \Big \rangle^2. 
\end{align*} 
\end{proposition} 

\textbf{Proof.} 
By Proposition \ref{second.derivatives.2}, we have 
\[\sum_{i=1}^2 \frac{\partial^2 Z}{\partial x_i \, \partial y_i}(\bar{x},\bar{y}) = \sum_{i=1}^2 (\lambda_i - \Phi(\bar{x})) = -2 \, \Phi(\bar{x}).\] 
Moreover, we have 
\[\sum_{i=1}^2 \frac{\partial^2 Z}{\partial y_i^2}(\bar{x},\bar{y}) = 2 \, \Phi(\bar{x}) \, \langle F(\bar{x}),F(\bar{y}) \rangle - 2 \, \langle \nu(\bar{x}),F(\bar{y}) \rangle = 2 \, \Phi(\bar{x}).\] 
Using these identities in combination with Proposition \ref{second.derivatives.1}, we conclude that 
\begin{align*} 
&\sum_{i=1}^2 \frac{\partial^2 Z}{\partial x_i^2}(\bar{x},\bar{y}) + 2 \sum_{i=1}^2 \frac{\partial^2 Z}{\partial x_i \, \partial y_i}(\bar{x},\bar{y}) + \sum_{i=1}^2 \frac{\partial^2 Z}{\partial y_i^2}(\bar{x},\bar{y}) \\ 
&= \Big ( \Delta_\Sigma \Phi(\bar{x}) - \frac{|\nabla \Phi(\bar{x})|^2}{\Phi(\bar{x})} + (|A(\bar{x})|^2-2) \, \Phi(\bar{x}) \Big ) \, (1 - \langle F(\bar{x}),F(\bar{y}) \rangle) \\ 
&- \frac{2 \, \Phi(\bar{x})^2 - |A(\bar{x})|^2}{2 \, \Phi(\bar{x}) \, (1-\langle F(\bar{x}),F(\bar{y}) \rangle)} \sum_{i=1}^2 \Big \langle \frac{\partial F}{\partial x_i}(\bar{x}),F(\bar{y}) \Big \rangle^2. 
\end{align*} 
This completes the proof.

\section{Proof of Theorem \ref{main.theorem}}

In this section, we describe the proof of Theorem \ref{main.theorem}. We first derive a Simons-type identity for the function $\Psi(x) = \frac{1}{\sqrt{2}} \, |A(x)|$.

\begin{proposition} 
\label{simons.identity}
Suppose that $F: \Sigma \to S^3$ is an embedded minimal torus in $S^3$. Then the function $\Psi = \frac{1}{\sqrt{2}} \, |A|$ is strictly positive. Moreover, $\Psi$ satisfies the partial differential equation 
\[\Delta_\Sigma \Psi - \frac{|\nabla \Psi|^2}{\Psi} + (|A|^2 - 2) \, \Psi = 0.\] 
\end{proposition}

\textbf{Proof.} 
It follows from work of Lawson that a minimal torus in $S^3$ has no umbilical points (see \cite{Lawson2}, Proposition 1.5). Thus, the function $|A|$ is strictly positive everywhere. Using the Simons identity (cf. \cite{Simons}, Theorem 5.3.1), we obtain 
\[\Delta h_{ik} + (|A|^2-2) \, h_{ik} = 0,\]  
hence 
\[\Delta_\Sigma(|A|^2) - 2 \, |\nabla A|^2 + 2 \, (|A|^2-2) \, |A|^2 = 0.\] 
The Codazzi equations imply that $|\nabla A|^2 = 2 \, \big | \nabla |A| \big |^2$. Consequently, we have 
\[\Delta_\Sigma(|A|) - \frac{\big | \nabla |A| \big |^2}{|A|} + (|A|^2 - 2) \, |A| = 0,\] 
as claimed. \\

\begin{proposition} 
\label{clifford}
Suppose that $F: \Sigma \to S^3$ is an embedded minimal torus in $S^3$. If 
\[\sup_{x,y \in \Sigma, \, x \neq y} \frac{|\langle \nu(x),F(y) \rangle|}{\Psi(x) \, (1 - \langle F(x),F(y) \rangle)} \leq 1,\] 
then $F$ is congruent to the Clifford torus. 
\end{proposition}

\textbf{Proof.} 
By assumption, we have  
\[Z(x,y) = \Psi(x) \, (1 - \langle F(x),F(y) \rangle) + \langle \nu(x),F(y) \rangle \geq 0\] 
for all points $x,y \in \Sigma$. For simplicity, let us identify the surface $\Sigma$ with its image under the embedding $F$, so that $F(x) = x$. Let us fix an arbitrary point $\bar{x} \in \Sigma$. We can find an orthonormal basis $\{e_1,e_2\}$ of $T_{\bar{x}} \Sigma$ such that $h(e_1,e_1) = \Psi(\bar{x})$, $h(e_1,e_2) = 0$, and $h(e_2,e_2) = -\Psi(\bar{x})$. Let $\gamma(t)$ be a geodesic on $\Sigma$ such that $\gamma(0) = \bar{x}$ and $\gamma'(0) = e_1$. We define a function $f: \mathbb{R} \to \mathbb{R}$ by 
\[f(t) = Z(\bar{x},\gamma(t)) = \Psi(\bar{x}) \, (1 - \langle \bar{x},\gamma(t) \rangle) + \langle \nu(\bar{x}),\gamma(t) \rangle \geq 0.\] 
A straightforward calculation gives 
\[f'(t) = -\langle \Psi(\bar{x}) \, \bar{x} - \nu(\bar{x}),\gamma'(t) \rangle,\] 
\begin{align*} 
f''(t) &= \langle \Psi(\bar{x}) \, \bar{x} - \nu(\bar{x}),\gamma(t) \rangle \\ 
&+ h(\gamma'(t),\gamma'(t)) \, \langle \Psi(\bar{x}) \, \bar{x} - \nu(\bar{x}),\nu(\gamma(t)) \rangle, 
\end{align*}
and 
\begin{align*} 
f'''(t) 
&= \langle \Psi(\bar{x}) \, \bar{x} - \nu(\bar{x}),\gamma'(t) \rangle \\ 
&+ h(\gamma'(t),\gamma'(t)) \, \langle \Psi(\bar{x}) \, \bar{x} - \nu(\bar{x}),D_{\gamma'(t)} \nu \rangle \\ 
&+ (D_{\gamma'(t)}^\Sigma h)(\gamma'(t),\gamma'(t)) \, \langle \Psi(\bar{x}) \, \bar{x} - \nu(\bar{x}),\nu(\gamma(t)) \rangle. 
\end{align*} 
In particular, we have $f(0) = f'(0) = f''(0) = 0$. Since $f(t)$ is nonnegative, we conclude that $f'''(0) = 0$. From this, we deduce that $(D_{e_1}^\Sigma h)(e_1,e_1) = 0$. An analogous argument with $\{e_1,e_2,\nu\}$ replaced by $\{e_2,e_1,-\nu\}$ yields $(D_{e_2}^\Sigma h)(e_2,e_2) = 0$. Using these identities and the Codazzi equations, we conclude that the second fundamental form is parallel. In particular, the intrinsic Gaussian curvature of $\Sigma$ is constant. Consequently, the induced metric on $\Sigma$ is flat. On the other hand, Lawson \cite{Lawson1} proved that the Clifford torus is the only flat minimal torus in $S^3$. Putting these facts together, the assertion follows. \\

We now complete the proof of Theorem \ref{main.theorem}. Suppose that $F: \Sigma \to S^3$ is an embedded minimal torus in $S^3$, and let 
\begin{equation} 
\kappa = \sup_{x,y \in \Sigma, \, x \neq y} \frac{|\langle \nu(x),F(y) \rangle|}{\Psi(x) \, (1 - \langle F(x),F(y) \rangle)}. 
\end{equation}
If $\kappa \leq 1$, then Proposition \ref{clifford} implies that $F$ is congruent to the Clifford torus. Hence, it suffices to consider the case $\kappa > 1$. By replacing $\nu$ by $-\nu$ if necessary, we can arrange that 
\begin{equation} 
\label{kappa}
\kappa = \sup_{x,y \in \Sigma, \, x \neq y} \Big ( -\frac{\langle \nu(x),F(y) \rangle}{\Psi(x) \, (1 - \langle F(x),F(y) \rangle)} \Big ). 
\end{equation}
We now define $\Phi(x) = \kappa \, \Psi(x)$ and 
\[Z(x,y) = \kappa \, \Psi(x) \, (1 - \langle F(x),F(y) \rangle) + \langle \nu(x),F(y) \rangle\] 
for $x,y \in \Sigma$. It follows from (\ref{kappa}) that the function $Z(x,y)$ is nonnegative, and the set 
\[\Omega = \{\bar{x} \in \Sigma: \text{\rm there exists a point $\bar{y} \in \Sigma \setminus \{\bar{x}\}$ such that $Z(\bar{x},\bar{y}) = 0$}\}\] 
is non-empty. Moreover, using Propositions \ref{key.ingredient} and \ref{simons.identity}, we conclude that 
\begin{align} 
\label{critical.point}
&\sum_{i=1}^2 \frac{\partial^2 Z}{\partial x_i^2}(\bar{x},\bar{y}) + 2 \sum_{i=1}^2 \frac{\partial^2 Z}{\partial x_i \, \partial y_i}(\bar{x},\bar{y}) + \sum_{i=1}^2 \frac{\partial^2 Z}{\partial y_i^2}(\bar{x},\bar{y}) \notag \\ 
&= -\frac{\kappa^2-1}{\kappa} \, \frac{\Psi(\bar{x})}{1-\langle F(\bar{x}),F(\bar{y}) \rangle} \sum_{i=1}^2 \Big \langle \frac{\partial F}{\partial x_i}(\bar{x}),F(\bar{y}) \Big \rangle^2 
\end{align} 
for every pair of points $\bar{x} \neq \bar{y}$ with the property that $Z(\bar{x},\bar{y}) = \frac{\partial Z}{\partial x_i}(\bar{x},\bar{y}) = \frac{\partial Z}{\partial y_i}(\bar{x},\bar{y}) = 0$.

\begin{proposition} 
\label{gradient.of.Psi}
We have $\nabla \Psi(\bar{x}) = 0$ for all points $\bar{x} \in \Omega$. 
\end{proposition}

\textbf{Proof.} 
Let us consider an arbitrary point $\bar{x} \in \Omega$. By definition of $\Omega$, we can find a point $\bar{y} \in \Sigma \setminus \{\bar{x}\}$ such that $Z(\bar{x},\bar{y}) = 0$. Since the function $Z$ attains its global minimum at the point $(\bar{x},\bar{y})$, the identity (\ref{critical.point}) gives 
\begin{align*} 
0 &\leq \sum_{i=1}^2 \frac{\partial^2 Z}{\partial x_i^2}(\bar{x},\bar{y}) + 2 \sum_{i=1}^2 \frac{\partial^2 Z}{\partial x_i \, \partial y_i}(\bar{x},\bar{y}) + \sum_{i=1}^2 \frac{\partial^2 Z}{\partial y_i^2}(\bar{x},\bar{y}) \\ 
&= -\frac{\kappa^2-1}{\kappa} \, \frac{\Psi(\bar{x})}{1-\langle F(\bar{x}),F(\bar{y}) \rangle} \sum_{i=1}^2 \Big \langle \frac{\partial F}{\partial x_i}(\bar{x}),F(\bar{y}) \Big \rangle^2 \leq 0. 
\end{align*} 
Since $\kappa>1$, we conclude that  
\[\Big \langle \frac{\partial F}{\partial x_i}(\bar{x}),F(\bar{y}) \Big \rangle = 0\] 
for each $i$. Using the identity (\ref{first.order.condition.a}), we deduce that 
\[0 = \frac{\partial Z}{\partial x_i}(\bar{x},\bar{y}) = \kappa \, \frac{\partial \Psi}{\partial x_i}(\bar{x}) \, (1 - \langle F(\bar{x}),F(\bar{y}) \rangle)\] 
for each $i$. Therefore, $\nabla \Psi(\bar{x}) = 0$, as claimed. \\

\begin{proposition}
\label{open}
The set $\Omega$ is open.
\end{proposition}

\textbf{Proof.} 
Let us consider an arbitrary pair of points $\bar{x} \neq \bar{y}$, and let $(x_1,x_2)$ and $(y_1,y_2)$ denote geodesic normal coordinates around $\bar{x}$ and $\bar{y}$, respectively. As in Section \ref{calc}, we can arrange that $\big \langle w_1,\frac{\partial F}{\partial y_1}(\bar{y}) \big \rangle \geq 0$, $\big \langle w_1,\frac{\partial F}{\partial y_2}(\bar{y}) \big \rangle = 0$, and $\big \langle w_2,\frac{\partial F}{\partial y_2}(\bar{y}) \big \rangle \geq 0$, where $w_1$ and $w_2$ are defined by 
\[w_i = \frac{\partial F}{\partial x_i}(\bar{x}) - 2 \, \Big \langle \frac{\partial F}{\partial x_i}(\bar{x}),\frac{F(\bar{x})-F(\bar{y})}{|F(\bar{x})-F(\bar{y})|} \Big \rangle \, \frac{F(\bar{x})-F(\bar{y})}{|F(\bar{x})-F(\bar{y})|}.\] 
Adapting the proof of Lemma \ref{algebra}, we conclude that 
\[\sum_{i=1}^2 \Big | w_i - \frac{\partial F}{\partial y_i}(\bar{y}) \Big | \leq \Lambda(\bar{x},\bar{y}) \, \bigg ( Z(\bar{x},\bar{y}) + \sum_{i=1}^2 \Big | \frac{\partial Z}{\partial y_i}(\bar{x},\bar{y}) \Big | \bigg ),\] 
where $\Lambda(x,y)$ is a continuous function on the set $\{(x,y) \in \Sigma \times \Sigma: x \neq y\}$, which may be unbounded along the diagonal. We next adapt the proof of Proposition \ref{key.ingredient} to obtain an estimate of the form  
\begin{align} 
\label{pde}
&\sum_{i=1}^2 \frac{\partial^2 Z}{\partial x_i^2}(\bar{x},\bar{y}) + 2 \sum_{i=1}^2 \frac{\partial^2 Z}{\partial x_i \, \partial y_i}(\bar{x},\bar{y}) + \sum_{i=1}^2 \frac{\partial^2 Z}{\partial y_i^2}(\bar{x},\bar{y}) \notag \\ 
&\leq -\frac{\kappa^2-1}{\kappa} \, \frac{\Psi(\bar{x})}{1-\langle F(\bar{x}),F(\bar{y}) \rangle} \sum_{i=1}^2 \Big \langle \frac{\partial F}{\partial x_i}(\bar{x}),F(\bar{y}) \Big \rangle^2 \\ 
&+ \tilde{\Lambda}(\bar{x},\bar{y}) \, \bigg ( Z(\bar{x},\bar{y}) + \sum_{i=1}^2 \Big | \frac{\partial Z}{\partial x_i}(\bar{x},\bar{y}) \Big | + \sum_{i=1}^2 \Big | \frac{\partial Z}{\partial y_i}(\bar{x},\bar{y}) \Big | \bigg ), \notag
\end{align} 
where $\tilde{\Lambda}(x,y)$ is another continuous function on the set $\{(x,y) \in \Sigma \times \Sigma: x \neq y\}$, which may be unbounded along the diagonal. Using the inequality (\ref{pde}) and Bony's strict maximum principle for degenerate elliptic equations, we conclude that the set $\Omega$ is open (see \cite{Bony} or \cite{Brendle-book}, Corollary 9.7). \\

Since $\Omega$ is open, it follows from Proposition \ref{gradient.of.Psi} that $\Delta_\Sigma \Psi(\bar{x}) = 0$ for each point $\bar{x} \in \Omega$. Hence, Proposition \ref{simons.identity} implies that $\Psi(\bar{x}) = 1$ for each point $\bar{x} \in \Omega$. Using standard unique continuation theorems for solutions of elliptic partial differential equations (see e.g. \cite{Aronszajn}), we conclude that $\Psi(x) = 1$ for all $x \in \Sigma$. Consequently, the Gaussian curvature of $\Sigma$ vanishes identically. As above, it follows from a result of Lawson \cite{Lawson1} that $F$ is congruent to the Clifford torus.

\end{document}